\documentclass[12pt]{notices}
\usepackage{amssymb,cite}
\usepackage{amsmath,graphicx}
\title{Mikhail Aleksandrovich Shubin\\
December 19, 1944 -- May 13, 2020}
\begin{document}
\maketitle
\begin{center}
\includegraphics[scale=0.7]{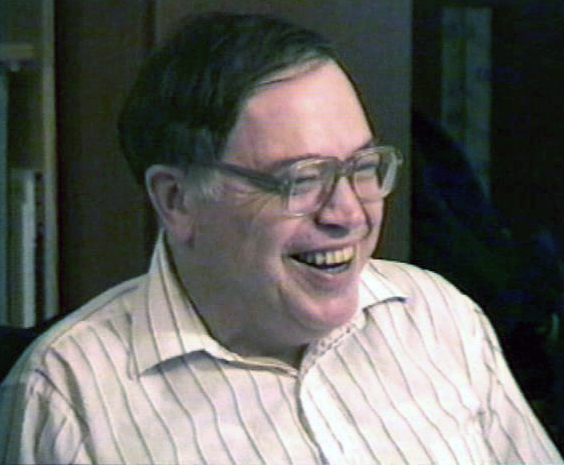}\\
\end{center}
Mikhail Aleksandrovich (Misha) Shubin, an outstanding mathematician, Fellow of AMS and Emeritus Distinguished Professor of Northeastern University passed away after a long illness on May 13, 2020. He was born on December 19, 1944 in Kujbyshev (now Samara) in Russia and brought up by his mother and grandmother. His mother, Maria Arkadievna, had graduated from the Department of Mechanics and Mathematics (mech-mat) of the Moscow State University (MSU) and became an engineer and later, after defending her PhD an Associate Professor at a Polytechnic Institute.

Initially Misha, who had an absolute pitch, was mostly interested and planned his career in music. In high school, however, he got involved with mathematics, competed successfully in olympiads, and then decided to become a mathematician. He was admitted to mech-mat of MSU in 1961 and became a pupil of a renowned expert in PDEs  M.~I.~Vishik. He defended in 1969 his PhD Thesis, devoted to the theory (including index) of matrix-valued Wiener-Hopf operators. In 1981 he defended his Doctor of Science degree (an analog of Habilitation), addressing the theory of operators with almost periodic coefficients, where he is one of the leaders.

Since then, M. Shubin had ventured successfully into a variety of areas of mathematics , demonstrating a remarkable width of interests and expertise. In his papers and books (totaling to around 140), he has made major and influential contributions to a variety of areas, including (but not limited to) operator theory, spectral theory of differential and pseudodifferential operators (especially operators with almost periodic and random coefficients), theory and applications of pseudodifferential operators and their discrete analogs, microlocal analysis, geometric analysis and analysis on manifolds (e.g., spectral theory on non-compact manifolds and Lie groups, index theory, general Riemann-Roch theorems, and invariants of manifolds), integrable systems, non-standard analysis, etc. Besides his many solo publications, he has co-authored works with his students, younger colleagues, as well as prominent experts, such as M.~Gromov, V.~Kondratiev,  V.~Maz'ya, S.~P.~Novikov, S.~Sunada, and others. For more details, one can consult with the memorial article \cite{UMN}.

Besides manifold major research papers, Misha has also co-authored several important survey articles. His book ``Pseudo-differential Operators and Spectral Theory'' \cite{PseudoBook}, first published in Russian in 1978, went through several English editions and still remains a popular source for studying microlocal analysis. Another remarkable book ``Schr\"odinger Equation'' \cite{Berezin}, which he co-authored with F.~Berezin, is also a well known and in some ways unique source. The AMS has just published his textbook ``Invitation to Partial Differential Equations'' \cite{PDEBook}, which is a revised and extended (edited by his colleagues M.~Braverman, R.McOwen, and P.~Topalov) version of his lectures given at the Moscow State University \cite{UrMatFiz}. He spearheaded and edited Russian translations  of the fundamental books by F.~Treves and L.~H\"ormander.

Misha had been participating in the famous I.~M.~Gel'fand's seminar from 1964 till his move to the USA in the beginning of 1990s. His notes for 25 years of the talks at the seminar have been, with the financial support from the Clay Institute, made available on the internet \cite{GelfSeminar}.

Misha Shubin was a remarkable lecturer and teacher, who has influenced lives and careers of many young mathematicians (from high school to graduate school and beyond). He was broadly educated, cheerful, and warmhearted person. In the 90s and 80s he helped several Soviet mathematician ``otkazniks.'' He also fought against notorious discriminatory admission practices at the Department of Mathematics of the MSU.

The memory of Misha Shubin, a mathematician, teacher, and a wonderful human being, will stay in the hearts of his colleagues, students, and friends.
\begin{input}

\begin{center}
 Maxim Braverman
\end{center}

I first met Misha Shubin when I was a graduate student in Tel-Aviv. Misha visited Israel and I asked him, as I understand now, a very basic question. Misha answered the question and gave me some references. To my surprise, the next morning he brought me several pages of neatly written hand notes with a review of the subject which even contained some complete proofs! ``I looked through the references I gave you," - said Misha, - ``And realized that they don't explain what you need well enough. So I wrote it down."

A few days later when we went to  the old city of Jerusalem I learned that Misha's  thoroughness  goes far  beyond his professional life.  I often took our guests to such tours and considered myself quite a decent amateur guide, but it turned out that Misha had done his homework. The guide-tourist relationship was almost completely reversed and my future tourist friends benefited greatly from the information I learned on this trip.

Several years later we were working together on a paper about self-adjointness of Schr\"odinger operators. At some point I wrote a 12 pages long draft and presented it to Misha. He seemed to like it but suggested working on some additions and generalizations. This happened again and again. Whenever I suggested to wrap things up and publish the paper, Misha answered: ``but it would be such a service to the mathematical community to analyze  this case as well". It was not just a phrase, he really cared about making a difference in the community. Half a year later we released a 50 page long manuscript which became my most cited paper. Most people who use it, cite one of these extra topics which we added ``to serve the community".

Misha loved books. He had a huge library. The walls of his apartment were covered with bookshelves. Once I asked him, how could he find a book among all these shelves.  He suggested that I choose any book in his library and ask him to find it. Each time I tried, he immediately went to the right shelf, took out the book, and told me when and why he bought the book, what is significant about it and why it stands in this particular shelf.

Misha always took very careful notes of every talk he attended and any mathematical discussion he had. His notes from the Gel'fand seminar in Moscow are well known. But there were many more. All of them were stored in file-cabinets in his office. Many times when we discussed math, he went to these cabinets and produced his notes from some talk he attended 10 years earlier where some of the questions we were discussing were addressed. Sometimes those were the notes from my own talks about which I have completely forgotten!

Over time, I learned that math was not enough to satisfy his curiosity. He wrote a very interesting essay about his visit to Egypt and conducted  nearly professional studies of whether Richard III really killed his nephews and why the Titanic drowned. It was always a pleasure to hear him talking about music, linguistics, history, and cultures of different places he visited.

Misha was always eager to help his colleagues, especially the  younger ones. His role in my own career is impossible to overestimate. When I already became an established researcher I mostly worked on the problems not directly related to Misha's mathematical interests. Still I often talked to him about my research. Only when Misha retired and I could not talk to him anymore, I realized how important those conversations were for my research and how much help he was ``invisibly" giving to me. I know many more mathematicians who owe a lot to Misha. The help he was giving to others throughout his career is a legacy in the field of mathematics that rivals his own research.

\end{input} 
\begin{input}
                                   \begin{center}
                                    Days with Misha\\
                                    Arnold Dikansky
                                   \end{center}

I have known Misha for more than 50 years. His mathematical knowledge was enormous and dedication to his profession -- legendary. Misha's advice to me when I was working on my PhD dissertation was crucial. Several distinguished mathematicians describe here his mathematical achievements and how he influenced his students and many other researchers. So, I will speak here only about the cultural dimension of Misha's multidimensional personality and about his family, M\&M (Misha and Masha).

We shared with Misha love for classical music and especially for Richard Wagner's operas. We had attended for many years various Ring operas at the Metropolitan Opera in New York. With the same level of thoughtfulness as he showed in his professional life, Misha approached upcoming performances, diligently studying librettos. When a new production of the Ring came out in 2012, I invited M\&M to come to New York, stay with me, and go to the performance of Die Walk\"ure. It was his last visit to New York, which he loved and especially appreciated its multicultural character.

It is interesting that many mathematicians enjoy deeply the classical music. Maybe the reason is that both music and mathematics touch something that lie very deep in human nature.
Both demand an enormous amount of time and effort in order to succeed. In his early life, when choosing the profession, Misha was torn between mathematics and music. He has achieved a professional summit in mathematics and enjoyed music and art in general as an amateur. Besides music, reading was Misha's another love. He had read the complete Shakespeare (in English, as well as in Russian translations).

I would like to say a few words about Misha's wife Masha. I have visited them on many occasions in their apartments in Moscow and Boston as well as, in the last several years, in the nursing home. The care which she gave to him was beyond what it seemed to be possible. It has lasted through many ups and downs in his condition, required tremendous perseverance, the control over her own emotions and reactions, and just physical strength.

Arnold Dikansky
Professor Emeritus, St. John's University New York

\end{input} 

\begin{center}
\textbf{Memories of Misha}\\
Leonid Friedlander
\end{center}


Probably I met Misha first time in 1975. I was still an undergraduate, and I was attending Misha's lectures on microlocal analysis. I remember him going in small details through all the calculus, and then he spent some time on the index theory. At first, his style looked too pedantic to me; only later I appreciated how important the details are.  That was one of the classes, on the basis of which he wrote his book "Pseudodifferential Operators and Spectral Theory", which is still one of the major textbooks on microlocal analysis.
A bit later, I started going to the Spectral Theory seminar that he was running in the Moscow State University. I spoke there several times. Quite often, Misha would ask a question that looked obvious. My first reaction then was: every idiot knows that; and then, in a couple of moments: I have no idea.

In 1979 I applied for emigration, and I was finally allowed to leave the Soviet Union in 1987. That 8 years of my life were not easy; most of the time I did not have a regular job. What made these years not only bearable but actually quite good was an enormous moral support that I received from my colleagues. In Misha's case, the support was not only moral. He arranged for me to translate a book into Russian; he himself edited the translation. Helping a political undesirable was not a neutral act. Misha could have big problems as a result of that, and he knew it well. As I learned later, I was not the only person who got a helping hand from him.  The honorarium that I received was substantial enough to solve my financial problems. Actually, when I was finally leaving, I paid the compulsory renouncing-of-citizenship fee from that money. Translating that book also gave me an opportunity to work closely with Misha. I was coming to Misha's apartment regularly, and we would go together through the new portion of the translation (there were no computers, no e-mail, everything was handwritten). Misha's motto was that the translation must me everywhere correct, regardless of whether or not the original is. If there is a mistake, it must be corrected, if there a gap, it should be filled. After two or three hours of work, usually at 10 or 11 p.m., we would go to the kitchen, and Masha, Misha's wife, would put khachapuries (Georgian bread filled with cheese) into the oven, and sometimes one or two Misha's friends who lived nearby would come, and, with tea and khachapuries, we would be talking till after midnight about fiction, politics, linguistics,... Misha had unlimited curiosity. That was true in mathematics, that was true in everything else.
\begin{center}
  \includegraphics[scale=0.8]{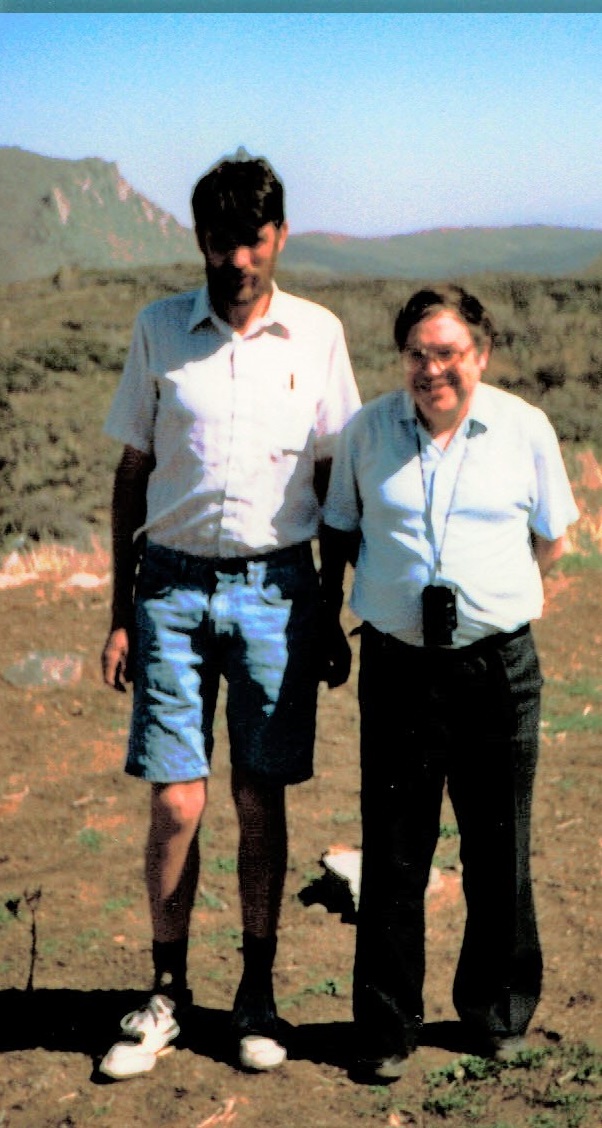}\\
  L. Friedlander and M. Shubin
\end{center}
Several years later, in May of 1991, Misha came to visit UCLA, where I was a post-doc at that time. We went together to San Diego.
On the way back, after having visited  Anza Borrego, rather late at night, I was driving on an almost empty freeway, with my wife sitting next to me, and Misha on the back seat. For half an hour or so, Misha was absolutely silent; we thought that he was sleeping. Then, he said: I have been looking at the road and thinking: how do they manage to maintain it so that all the bulbs are working. What bulbs? Turned out that he thought that reflectors are electric bulbs (there were no reflectors on roads in Russia, so he saw them first time). What he saw did not make sense, and he was trying to solve the problem.

Misha had a strong moral compass. He believed that certain things are worth fighting for. Back in the USSR, he was not a dissident, but showing  independence was also fighting the system. As a result of that, he had serious career problems. On the other hand, he had tolerance to human weaknesses. I remember well a conversation that we had in his kitchen in Moscow. We were talking about one specific case, and I was a hardliner, as most young people are. Misha was trying to convince me that, in the particular situation, one should show more understanding and be more forgiving. We respectfully disagreed then. Now, after almost 40 years, I think that he was probably right.

I met Misha many times and in different places. When I was on Sabbatical in MIT in 2004--2005, we talked a lot. He invited me to give a talk at the Basic Notions Seminar that he was running in the Northeastern University. That was his wonderful invention. I spent much more time preparing that talk than I would usually spend for a research talk. Last time I saw him in 2009 at the conference for his 65th birthday. I was fortunate to have known Misha and to be his friend.


\begin{input}
\begin{center}
\textbf{A few words about Misha Shubin}\\
Misha Gromov
\end{center}

What was striking in Misha were his unlimited fearlessness and his insatiable
love for knowledge. His eyes were radiating with light when he was coming across
a new idea, especially a mathematical one.

His knowledge and his understanding of mathematics was extraordinary, his
sensitivity to everything in the world was astonishing and he generously shared
his knowledge and his ideas with everybody.

I cherish the memory of spending time in his company,
I wish there were more such people as Misha among us.
\end{input}

\begin{input}
\begin{center}

\textbf{My friend Misha Shubin}\\
Victor Ivrii
\end{center}

I met Misha Shubin  for the first time in early 70s when I, a graduate student of the Novosibirsk State University, came to the  Moscow University to give several talks.  Misha had already defended his dissertation;  we instantly became friends. At that time ``Vishik's seminar'' was the main center  for Microlocal Analysis in USSR, and Misha was a very prominent member of this seminar.

I have visited Moscow many times after that, from Novosibirsk and later from Magnitogorsk, often stayed in Misha's  apartment; we met at several conferences in the Soviet Union and abroad. In particular, in October of 1981 we went together to the conference ``Differential Geometry and Global Analysis'' (Garwitz, GDR=East Germany). We were part of the Soviet delegation  (in the USSR that was the way to send scientists abroad). After many rejections,  it was my first visit abroad, and one of  Misha's  first such trips, after many rejections as well. At that time, a simple conference trip abroad required a very long chain of permissions, and  majority of academics never ventured abroad.
And finally we got one!  We felt like inmates let from our cell to the inner yard. It was a kind of a test: were we ``politically mature'' and could we be trusted to visit countries which were not the part of the Eastern Block? We failed this test miserably and due to a report written by two members of the same delegation (one of whom was an official snitch) that would have been our last trip beyond borders of the USSR, if not the  Perestroika.

That was a large international conference, attended by a number of our colleagues from Western Europe. We had known one another scientifically, but never met in-person  (one should remember that there was no email available for us, leave alone Skype or Zoom). .No wonder that, unlike some more senior members of our group who kept their own company and spoke exclusively Russian, we talked a lot and not only about mathematics. We spoke  English, not understood by  our ``colleague,'' who asked ``What are you talking about?'' We responded ``It is of no interest to you, it is mathematics''.
%

In the early 1990s we were roommates  at the Institut des Hautes \'Etudes Scientifiques. For both of us it was the first visit beyond Eastern Block (which already had collapsed). Everything was new and very exciting, from streets of Paris to museums. Misha's biggest complain was ``We will leave before we taste one tenth the yummy staff in the grocery stores!'' Later I visited him and his family several times in Boston and stayed at his place. He also was visiting me and my family  in Toronto, and also stayed in our house.

Misha has played a very important role in my mathematical life: it was he and Boris Levitan who suggested that I try to prove Weyl's conjecture in 1978.

Misha was a remarkable mathematician, very broad, with a great mathematical taste.  He was collaborating with many mathematicians on the large variety of topics.
His  book ``Pseudodifferential operators and spectral theory'' played very important role in promoting Microlocal Analysis and  related topics in Soviet Union and despite availability of other books it has become one of my favorites,   as well as his and Felix Berezin's book ``Schr\"odinger Equation''.  He also was an editor of several volumes of the series ``Encyclopaedia of Mathematical Sciences. Modern Directions".

Due to space limitations, I will mention just a couple of his many results. He and Vladimir Tulovsky have developed the \emph{Method of the almost spectral projector} in the theory of spectral asymptotics, which was in some sense an intermediate between Weyl's and Courant's variational approach and the Tauberian approach due to T.~Carleman. In this method a pseudodifferential operator was directly constructed that was close to the spectral projector.

In his paper with D.~Schenk ``Asymptotic expansion of the state density and the spectral function of a Hill operator" the \emph{complete spectral asymptotics} was derived for the integrated density of states for one-dimensional Schr\"odinger operator with periodic potential. This result was much later generalized to multidimensional case by L.~Parnovski and R.~Shterenberg.  The most amazing thing there was that the complete spectral asymptotics does not contradict to the fact that  for one-dimensional  operators with generic potentials all spectral gaps are open (but due to complete asymptotics these gaps are very narrow and  asymptotics itself  cannot not be differentiated).

\end{input}

\begin{input}
\begin{center}
\textbf{Yuri Kordyukov}
\end{center}

I first met Misha Shubin back in 1980, when I was a sophomore at the Department of Mechanics and Mathematics of the Moscow State University (MSU). At that time he was, of course, not Misha for me, but Mikhail Aleksandrovich. He became Misha much later,  sometime in the late 1990s. Shubin lead our recitations on the ordinary differential equations course. It was the time to choose a research advisor. I had been thinking for a long time, attended lectures of various mathematicians, and in the end my choice settled on Shubin.

As far as I understand now, at that time Shubin had just defended his Doctor of Sciences dissertation  and was looking for new directions of research. In 1981, he wrote a paper on the spectrum distribution function for transversally elliptic operators on a compact manifold equipped with a Lie group action. In 1983, Shubin wrote a paper with his graduate student Meladze devoted to pseudodifferential calculus on Lie groups. At that time he had a group of students working on analysis on noncompact manifolds and index theory. For instance, Alexander Efremov studied index theory on coverings of compact manifolds, Dmitry Efremov - spectral theory on hyperbolic plane. Brenner worked on analysis on manifolds with boundary. The weekly seminar led by Shubin met on Fridays. We studied the classic works by Atiyah-Bott-Patodi, the newly emerging works of Bismut, the book by Guillemin and Sternberg on geometric asymptotics, etc. I remember Shubin giving a course on the Atiyah-Singer index theorem. But his research interests were much broader. At the same time, in the first half of the 1980s, he published his work with Bondareva about solutions of the Cauchy problem for nonlinear integrable equations, a survey with Zvonkin on nonstandard analysis, papers with Schenk on complete asymptotic expansions of the density of states for a periodic Hill operator, and many others. He also edited the Russian translations of the two-volume Treves’ monograph, the four-volume monograph by H\"ormander, and the book by Rempel and Schulze on the index theory for elliptic boundary value problems. I have heard that after writing his book on pseudodifferential operators, he began writing a book on Fourier integral operators, and had already written 150 pages. But when he learned about the publication of Treves’ book, he stopped his own writing on the subject and with great enthusiasm switched to translating the Treves’ book. Shubin had always been very busy and worked hard. It was his natural state. He expected the same from me.

As a subject for my research, Shubin suggested to study analogs of his results on transversally elliptic operators for the simplest example of a noncompact Lie group --- the group of real numbers. After a while, during our discussions, the idea came up to look at the 1979 paper on noncommutative integration theory by Connes, and I switched to studying it. Later, manifolds of bounded geometry came into my view. The results on uniformly elliptic differential operators on manifolds of bounded geometry formed the basis of the first chapter of my Ph.D. thesis, defended under Shubin's supervision in 1988. The second chapter was devoted to tangentially elliptic operators on foliated manifolds. In 1988, after completing my graduate studies, I left for Ufa, where I found a position at the Aviation Institute. Shubin was not sure whether I would continue to work with him, and thus recommended to me some mathematicians with whom I could work in Ufa. In the end, I decided to continue the research that I had started under Misha's direction. It was already the beginning of the 1990s, the period of big changes in Russia. In 1991, Shubin went to MIT for a year and then found a permanent position at the Northeastern University.

After leaving for the United States, Shubin would rarely come to Moscow. I had met him several times in his Moscow apartment in the 1990s. In 1998, our collaboration resumed to continue his work with Gromov and Henkin on $L ^2$-holomorphic functions on coverings. Regretfully, this project has never been finished. In 2001, Shubin found an opportunity to invite me to Boston for a month. This was my first acquaintance with the United States, and Misha was happy to introduce me to various nuances and peculiarities of American life. In Boston, he invited me to join his project with Mathai to study the asymptotic behavior of the spectrum of a periodic magnetic Schr\"odinger operator in semi-classical limit. During my stay at Boston, within the framework of this project we managed to solve a problem that had applications to the quantum Hall effect. These results were published in 2005 in the Crelle’s Journal. It was our only joint paper with Misha. This research triggered my continuing interest in the magnetic Schr\"odinger operators. Then we parted our ways again. Shubin was very involved in his collaboration with Kondratiev and Maz'ya on the spectral theory of the magnetic Schr\"odinger operator. When I was back to Boston in 2001, he showed me some parts of these works that involved very subtle analytical arguments. I last met Shubin in 2009 in Boston at the Conference in honor of his 65th birthday. At that time he was already seriously ill, and we did not talk much about mathematics.

Shubin had a great influence on my life, taught me a lot. He has largely determined my path in mathematics. He shared with me his vision of mathematics, value criteria of mathematical results, his principles of joint scientific research and interaction with other mathematicians. I remember that in my second year at MSU, when I was just choosing my path in mathematics, he explained to me the unity of mathematics. He said that mathematics is one big interconnected area. Therefore, it does not matter in which part of mathematics to begin your scientific activity.  In the future it will be possible, if necessary, to move from one part to another, with the exception of, perhaps, only a few isolated islands, and the main goal is not to get stuck on one of them. He had a huge library of mathematics books (about 5 thousand, as far as I remember), which he had collected back in the Soviet times. He has managed to move it to America. When I was in Boston in 2001, I stayed at Misha's place and had accommodation in his library (which was also his study). I lived for a month in the midst of all these wonderful books. He told me; ``Yura! Why do you need to go to the university? You have everything here for work - a computer, books, the Internet.'' He was in the habit of working at home, developed through long years of work at Moscow State University. MSU professors did not have individual offices, and thus they usually worked from home. He taught me the skill of making presentations (``you should write on the board, starting in the upper left corner'' and so on), explained the cut and paste techniques for writing mathematical texts in that distant pre-computer time. He also helped me with various life tasks (be it job search after graduation or choosing a typewriter). In his turn, he had mentioned many times that he had learned a lot in mathematics from me. He has never raised his voice, always tried to be extremely tactful and delicate, not to impose his opinion. He was principled, honest, demanding of himself, very humble and sincerely devoted to mathematics. So he will remain in the memory of many people who knew him.

  \end{input}

\begin{input}
\begin{center}
\texttt{\textbf{Misha Shubin, dear friend and colleague}}\\
Peter Kuchment
\end{center}
It is hard to believe that Misha Shubin, such a vibrant, brilliant, and friendly person, is not with us anymore.

We have met probably somewhat more than a half-century ago, when I was still an undergraduate, and he - a graduate student. We met at one of the first famous in the former Soviet Union Voronezh Winter Mathematical Schools.
These were unimaginable feasts of mathematics, involving hundreds mathematicians (from undergraduate students to the leading experts) from the whole Soviet union. The discussions among younger participants like us were very intensive, every lecture being taken apart. This is where many friendships like ours with Misha and collaborations started. Ours has lasted throughout half of a century and different cities and countries. Back in Russia, I visited Moscow and Misha visited Voronezh with seminar talks and discussions. We have never done any joint work, but our mathematical paths were often closely parallel. For instance, at about same time when Misha worked on almost-periodic PDEs, I did the periodic ones; then this repeated with us both looking at operators on homogeneous spaces, discrete problems, etc. His work has influenced most of mine. His graduate school days paper on holomorphic projections on holomorphic sub-bundles of a trivial Banach bundle, which used intricate and recent at that time results of SCV going back to the work of H.~Grauert, has been used by me and many other operator theorists since. I have learned a lot from his work on almost-periodic PDEs. In 1975, being on sabbatical I attended his wonderful lectures on microlocal analysis at Moscow State University. Later on, they turned into his famous book \cite{PseudoBook}. The lectures were so meticulously prepared that it is no wonder they could be turned into a book rather easily. This was my practically first significant exposure to microlocal analysis, pseudodifferential operators and such, which turned out to be crucial for the work I started doing a couple of years later. I still rip the benefits of his lectures in most areas of my interest, especially in medical imaging and PDEs. Very recently I worked with my student on extending further Gromov's and Shubin's Riemann-Roch type results.
Moreover, Shubin's influence on my growth as mathematician and my research extended much further than direct links that I have just indicated. Our discussions, his lectures and his publications have taught me a lot (and I do not think I was his best student).

When I and my family were leaving the fSU in 1989 (see the photo of the farewell party),
\begin{center}
  \includegraphics[scale=0.6]{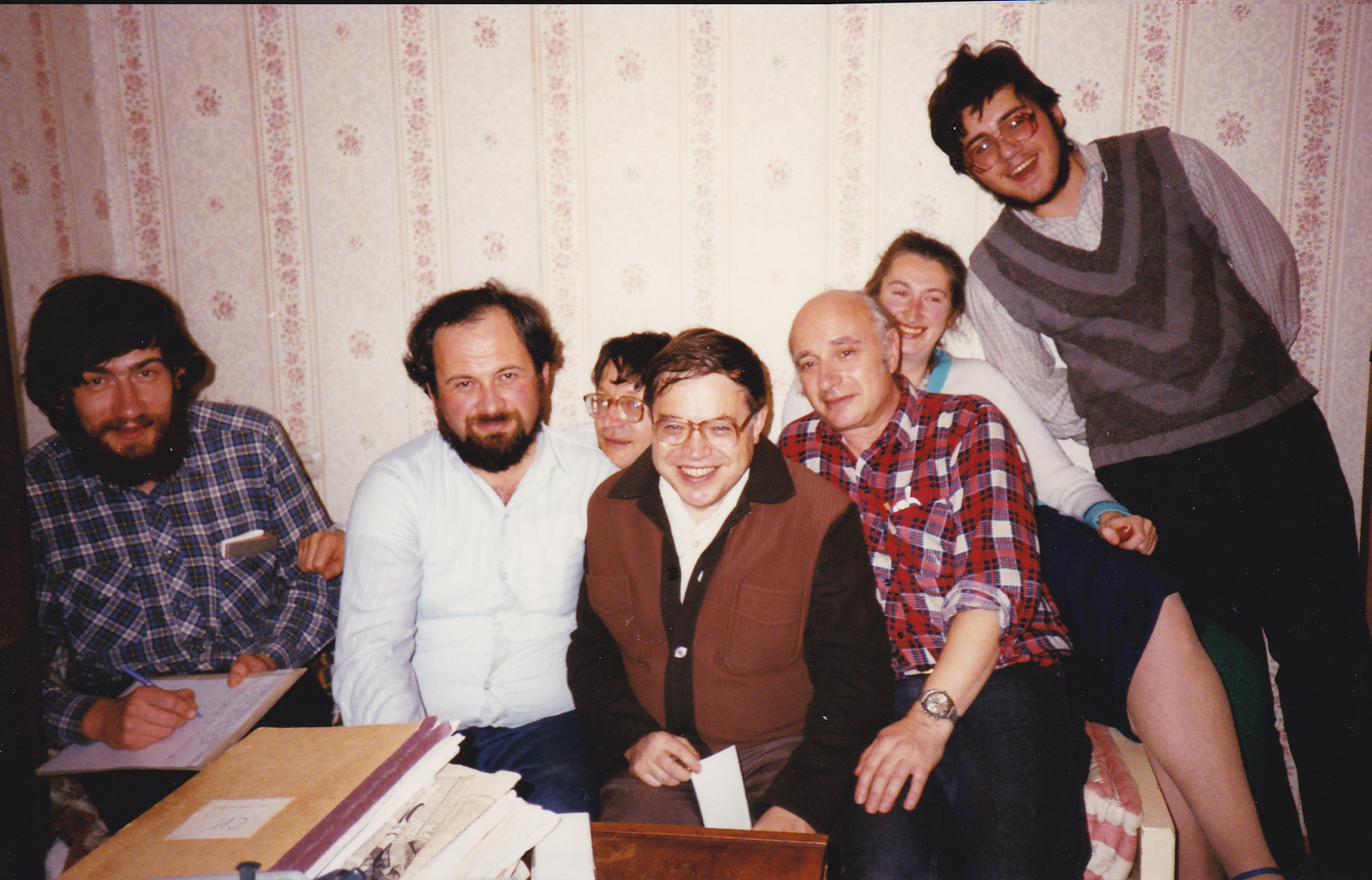}\\
  S.~Orevkov, P.~Kuchment, V.~Lomonosov, M.~Shubin, V.~Lin, E.~Landis, L.~Posicelskii
\end{center}
 we all thought that we were parting forever, with no chance to meet ever again. Fortunately, the times have changed, and Misha also moved to the USA. We kept communicating a lot and visiting each other, me coming to Boston, and him visiting Wichita, KS and College Station, TX, where I have worked.

Misha has played a very significant role in my life in many other regards. He was an extremely broadly educated person and was happy to share his knowledge. I have read quite a few fiction and non-fiction books at his recommendation. I had not studied English before 1985 at all (my foreign language was German). In 1985, during my another sabbatical in Moscow, Misha recommended me a wonderful immersion English class, which he had taken previously. He also made sure I could get into it (it required a recommendation). This has played a major role in my future life, especially after my emigration to the USA.

In difficult and corrupt Soviet times, he was one of the staples of moral fortitude, which often did not make his life any easier, to say the least. He was always ready to  help people in difficult circumstances. I remember once, when I lived in Kansas and he was on a business trip to Mexico, he called me from there to give an advice about universities choice for my daughter, who was graduating from the high school at that time.

Misha Shubin was a wonderful person in all regards: as a mathematician, teacher, writer, and most importantly, as human being. We all will miss him thoroughly. We are grateful to the fortune for letting us know him.
\end{input}

\begin{input}

\begin{center}
\textbf{Misha Shubin, An Exceptional Erudite}\\
Vladimir Maz'ya
\end{center}

It was by pure chance that Misha Shubin and I happened to share a compartment in a night train from Moscow to Ruse (Bulgaria) in 1975. We both were participants at a conference and for both of us it was the first trip abroad. As novices in traveling  outside the USSR, we could only be allowed a trip to a country within the Soviet bloc.

The train trip lasted almost three days and we had plenty of time to speak about mathematics as well as about our life in two different mathematical communities. He was telling me about the Division of Differential Equations at Moscow University, where he worked, and I told him about my research position at the Research Institute within the Math Department of Leningrad University, at the laboratory headed by Solomon G. Mikhlin. I reminded Misha that I had been present at two of his talks at Vishik's seminar in Moscow. Seven years younger than me, Misha  surprised me by his broad knowledge in mathematics.

As it often happens, our acquaintance was warmed by shared food. There was no restaurant on the train and we offered each other "delicacies" like cooked eggs, boiled chicken, tomatoes and cucumbers. The only thing provided by the train hostess was hot tea in  mornings and evenings. It is interesting that I remember this train trip with Misha more vividly than the conference we attended.

Later on, we would meet regularly at the Petrovsky winter conferences in Moscow, at which I regularly participated. However, a close contact became possible only much later, when we both emigrated.

In 2003, Bob McOwen invited me to spend six months as a visiting professor at the Northeastern University in Boston. My wife Tanya had a lot of teaching duties, while I was privileged to give just a special topic course for graduate students for two hours a week and participate in the seminar, where I gave several talks.

Besides doing research with Bob McOwen, who was the Chairman at that time, very soon I became involved in the following interesting mathematical problem.
Let $ -\Delta + V$ be the Schr\"odinger  operator in $L^2(\Bbb{R}^n)$, where $V\in L^1_{loc} (\Bbb{R}^n)$ is bounded from below; $n\geq 2$.  A.M. Molchanov discovered in 1953  that for some constant  $\gamma>0$ the following  condition  is  equivalent to the discreteness of the spectrum: for every $d>0$
\begin{equation}\label{1}
\inf_F \int_{Q_d \backslash F} V(x) dx \to \infty, \qquad Q_d \to \infty,
\end{equation}
where $Q_d$ are open  cubes with edge length $d$, and the compact sets $F\subset \overline{Q_d}$ are negligible in the sense that $cap(F) \leq \gamma\, cap(Q_d)$.  I.M. Gel'fand had raised a question about the best possible constant $\gamma$.  Misha and I answered this question by proving that $\gamma$ can be taken arbitrarily in $(0,1)$. Moreover, we showed that the negligibility condition can be weakened to
\begin{equation}\label{2}
cap(F) \leq \gamma(d)\, cap(Q_d)
\end{equation}
and described all admissible functions $\gamma: (0, \infty) \to (0,1)$. They must satisfy
\begin{equation}\label{3}
\mathop{\hbox{\rm lim\, sup}}_{d\to 0}\, d^{-2} \gamma(d) = \infty.
\end{equation}

All conditions (1) with functions $\gamma$ satisfying (2)-(3) are necessary and sufficient for the discreteness of the spectrum of the  Schr\"odinger operator. Our solution of this fascinating problem appeared in  "Annals of Mathematics" in 2005.
\begin{center}
  \includegraphics[scale=0.5]{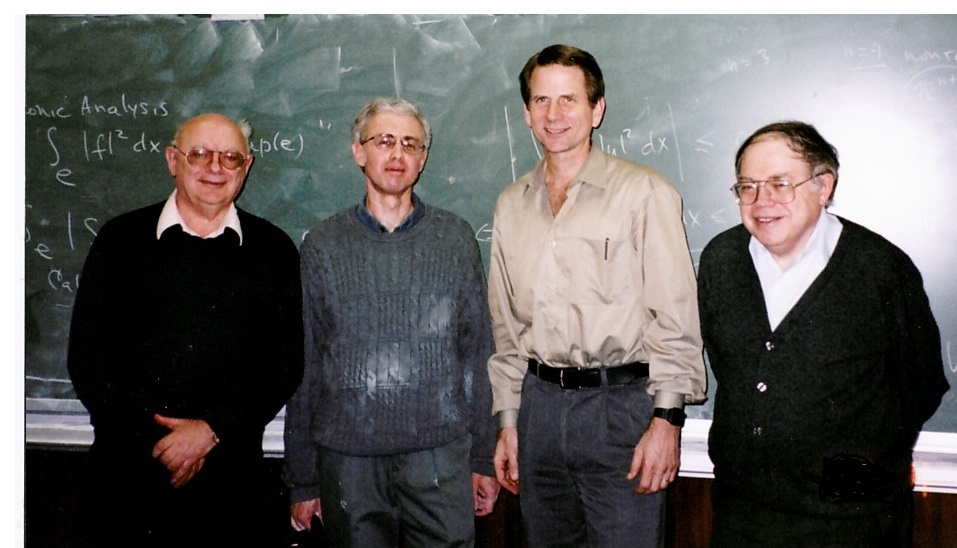}\\
  \textbf{Left to right: V. Maz'ya, I. Verbitskii, R. McOwen, M. Shubin}
\end{center}
Another question we succeeded in answering  became the topic of our paper ``Can one see the fundamental frequency of a drum?'' (Letters in Mathematical Physics, 2005). The main theorem here states that the bottom of the Dirichlet Laplace spectrum $\lambda(\Omega)$ for an open set $\Omega\subset (\Bbb{R}^n)$ is equivalent to $r^{-2}_\Omega$, where $r_\Omega$ is the so-called interior capacity radius of $\Omega$. The reason for the resemblance of the article's title to the famous Marc Kac's problem about isospectral drums can be explained as follows. If an eye could filter out sets of ``small'' capacity, then one could recover ``by looking'' the lowest Dirichlet Laplacian's eigenvalue, or at least get its reasonably good estimate.

When Volodya Kondratiev, who worked at the same division at Moscow State University where Misha used to work, came for a short visit to Misha in Boston, three of us became involved in a  problem on spectral properties of the magnetic Schr\"odinger operator. The first publication by the three of us appeared a year later in ``Communications in PDEs.'' The second one was finished later and published in 2009.

We worked with Misha intensively for six months, meeting both at the Department of Mathematics and at his place on weekends. My wife Tanya and I rented an apartment within five minutes' walk from Misha's home. The proximity to Misha's place was the only advantage of our rather modest apartment, which was made tolerable with the generous help from Masha Shubina (Misha's wife), who offered kitchen utensils and other necessities.

The six months in Boston were fruitful. Together with Misha we wrote four articles. However, when we left Boston, the distance and the lack of technical means of communication like Skype, made our further collaboration more difficult. We spent a lot of time on the phone and exchanging emails, but our contacts somehow faded away. When we discussed possible directions of further collaboration with the starting point being our ``Annals'' paper of 2005, it became clear that we were rather different in our preferences. Misha, being an exceptional erudite, wanted to widen the scope of our efforts, while I usually prefer going after necessary and sufficient conditions in a problem. This difference in tastes was difficult to overcome. It is a pity that we have not had an opportunity to work together since then.

\end{input} 
\begin{input}

\begin{center}
\textbf{Memories of Misha}\\Robert McOwen \\
Department of Mathematics
\end{center}

I first met Misha Shubin in Spring of 1987 at a mathematics conference in Oberwolfach, Germany. Both of us were giving talks at the conference. I don’t remember what topic Misha talked about, but I remember him most clearly from one of the evening socials towards the end of the conference. It was his first time out of the Eastern Bloc, and he had the opportunity to meet in person many mathematicians whom he had only known by their work. I remember him getting very emotional and, with a big smile on his face, saying “I am so happy!”

It was five years later that Misha joined the Math Department at Northeastern University in Boston. I remember his arrival having a big impact on our research profile and our graduate program. In Fall 1993, he started the Basic Notions Seminar. Topics of seminar talks varied, but the rule that Misha enforced relentlessly was that all stated mathematical results must be proved from elementary principles. Speakers who gave a series of talks were subject to another rule that Misha also enforced: subsequent talks must review all results from previous talks before proceeding to new material. The third rule that Misha imposed was that questions were strongly encouraged and must be answered in full detail. Misha usually set the pace himself by asking lots of questions that slowed the speaker down. Needless to say, not a lot of material was covered in any given talk, but the audience always came away from the Basic Notions Seminar with a firm understanding of what they had heard.

In 1995, my colleague Christopher King and I co-organized with Misha a Special Session on ``PDEs in Geometry and Mathematical Physics'' for an  AMS Conference that was held at NU.  Misha had a huge list of names to suggest as potential speakers, and we had to select carefully: they were all eager to come to NU to give a talk in a Special Session that he had co-organized. Misha served as Graduate Director of the Math Department from 1995 to 1997. During this period of time, he introduced many new courses to our graduate curriculum, like a course on the Atiyah-Singer Index Theorem. I think that the course was only taught once, by Misha himself, but it shows the ambitious standard that he aspired to for our graduate program.

Misha became a Matthews Distinguished University Professor at Northeastern in 2001. He was the first Distinguished Professor appointed in the Mathematics Department and remained the only one until Andre Zelevinsky was given the honor after his death in 2013. Misha earned the honor by being a highly respected researcher and author of many published articles and books. I have always been impressed by the number of mathematicians who cited his 1978 book, {\it Pseudo-differential operators and spectral theory}, as having a big impact on their research. Consequently, I was honored when the AMS asked me to help bring out the textbook on partial differential equations that Misha had been working on for the last years of his productive life. I had not seen the book before, but as I looked through it, I was impressed not just with the content but with the style of writing: he provided rigorous proofs but written in a friendly way that made them easily understood. I felt it was important that  the book be published, so enlisted the assistance of my colleagues Maxim Braverman and Peter Topalov. I was glad that we were able to complete the editing process and have {\it Invitation to Partial Differential Equations} appear in the AMS list of publications slightly before Misha’s death this year. I hope it will serve as a lasting reminder of the gentle man who was an inspirational figure in mathematics.

\end{input} 
\begin{input}

\begin{center}
\textbf{Coincidence or Destiny\\
In memory of Misha}\\
Toshikazu Sunada
\end{center}

For Fortuna, the goddess of destiny in Roman religion, there exists no difference between ``coincidence" and ``destiny," but for human being who cannot judge the appearing spots of a die in advance, the encounter with people, for instance, is neither more nor less than a ``chance".

It was in 1993 when I met Misha for the first time, on the occasion of the Japan-U.S. Mathematics Institute (JAMI) Conference ``Zeta functions in number theory and geometric analysis" held at Johns Hopkins University. I already knew him by name as a leading figure in the field of analysis, especially, of elliptic PDEs pertinent to physics and geometry. One day, in the tea break of the conference, he approached me, and introduced himself with a gentle and  cordial way of talking. Although I cannot recall exactly the content of our conversation then, it is certain, however, that we discussed some mathematics and that his interest in part turned out to have something in common with mine. According to my poor memory, our common interest was the spectra of ``magnetic Schr\"{o}dinger operators" in both continuous and discrete cases. His approach was analytical in nature, while mine was geometrical, if anything. In any event, I was pleased to know that Misha had interest in my research, for I thought of what I was studying at that time as being isolated in the mathematical community in Japan. Indeed, our conversation gave me a great deal of motivations and encouragement. But, that is not all; I was attracted by his personality during the course of the conversation. He really is a caring and warm-hearted person.

The next chance to meet him came in 1996. We both participated in the symposium on ``Partial Differential Equations and Spectral Theory" held at Durham University, UK. Taking this opportunity, we discussed not only on mathematics, but also on non-mathematical subjects. I was very much impressed with his broad knowledge of humanities in general. Above all, he showed his keen interest in Japanese culture. He has read some Japanese contemporary literature, through which he got well acquainted with Japanese mentality and behaviour.

In 2005, I invited him to Japan to do a joint research. At that time, I was trying to make the classical $T^3$-law of specific heat rigorous by means of discrete geometric analysis. During his stay in Tokyo, where he enjoyed longed-for Japanese life, we could accomplish a joint work entitled ``Mathematical theory of lattice vibrations and specific heat" that was published in Pure and Appl. Math. Quaterly, a newly established journal. It is my great honour to be a collaborator of Misha.

The last chance to have met him was in 2009. I was invited to the Conference in Honour of the 65th Birthday of Mikhail Shubin, ``Spectral Theory of Geometric Analysis" held at Northeastern University, Boston. I could meet his beautiful family in the birthday party at his house.  He looked happy  and contented, surrounded by his family and many friends. Eleven years since then, I still remember his warm smiling face. I wanted to see him again and to have a pleasant conversation about a wide variety of things, but his illness precluded me to do it.  Last May, I was deeply saddened when  I heard that he passed away. I really bear a grudge against Fortuna for not giving me the chance to fulfill my desire.

To conclude this eulogy, I would like to add a few statements, in relation to his personality. On $11^{\rm th}$ March, 2011, Japan experienced Great East Japan Earthquake. The tsunami triggered by the earthquake killed more than 15000 people along the east coast of Tohoku region. Misha was worried about my family (in my memory, he gave me an international call from Boston). Fortunately, we were not affected by this disaster since Tokyo is far away from Tohoku. What I wish to say is that he not only confirmed my family's safety, but also made a donation to the Mathematics Institute of Tohoku University to assist a few teaching staff affected by this disaster (I had been a member of the institute for ten years). The members of the institute very much appreciated his kind offer and heart-warming message. I, too, am thankful to Misha for his solicitude.

Well Misha, my friend and collaborator, your sincerity, kindness, warmness, generosity, caring for others, smiles, humours will stay with me forever.  May you rest in eternal peace.
I extend my most sincere condolences to Misha's family.

\end{input}

\begin{input}
\begin{center}
\textbf{Misha Shubin, a few personal memories}\\
Alexander Zvonkin
\end{center}

Misha Shubin was the most important person in my life after my family -- my wife
and children.

At a certain point I was going through a difficult period: I had not been doing math
for almost ten years. As a rule, it is next to impossible to restart after such a
long interruption --- unless you get help from someone. It was a good fortune for me
that Misha and I became neighbors. It started with small things: while taking the same
bus quite by chance, we discussed various mathematical subjects. Then Misha proposed
that I translate a French mathematical paper. Then he lent me a book on nonstandard
analysis, the subject that interested him at that time. This way, step by step,
we finished by co-authoring quite a substantial paper. But that was not all:
Misha convinced me to start participating in the weekly seminar headed by I.~M.~Gelfand,
and my next paper was greatly influenced by the latter.

Later, while visiting the Institut des Hautes \'Etudes Scientifiques in Bures-sur-Yvette,
Misha recommended me as one of the possible guests, and that visit was absolutely crucial
for my life and hugely predetermined it.

I could easily continue with this list but then it would be about myself, not about Misha.

Misha was a very good person. I really do not know how to say it in a more elegant way.
He was just a very good person. I remember once talking with a dissident friend of mine
who claimed that in our country of omnipresent denunciations and the overwhelming power
of the KGB you could not trust a single person: anyone could turn out to be a snitch.
As a counterexample, I told him that Misha was a person whom one could totally trust.
He thought for a while and then, despite of his habit of always insisting on his point of view, said: ``I think you are right this time; it is impossible to imagine Misha as an informant.''

Misha saw that my way of writing mathematical papers was unprofessional, and he found a way to let me know it without hurting my feelings. He asked me whether I knew Asimov's stories about robots and then, proceeding from the famous Three Laws of Robotics, formulated similar laws of good math paper writing:
\begin{quote}
{\sc Three basic principles for writing a mathematical paper}
\begin{enumerate}
\item	The text of the paper should be mathematically correct.
\item	The text should be written as rapidly as possible except when this would conflict with the First Principle.
\item	One should try to make the text as perfect as possible (comprehensible, accessible, self-contained and complete) except when this would conflict with the first two Principles.
\end{enumerate}
\end{quote}
It was almost forty years ago, but I still keep the sheet of paper with these
principles written by him.

Our paths have diverged. Misha moved to Boston, while I moved to Bordeaux. Misha has visited us in France five times. Those were purely personal visits, though we could not avoid talking about mathematics. His casual remarks, his questions, or his surprise at certain phenomena helped me to better understand my own results.

He has never made an impression of being in a hurry, but worked amazingly fast.
But no matter how fast he worked, it was never at the expense of meticulousness
(remember the first principle: ``The text should be mathematically correct'').

When we discussed the translation of the above-mentioned French article, he asked me
how long I thought it would take. At that time I had a full-time job at an applied research institute, but hoped to be able to do the translation in a month.
``A month!'' -- Misha was clearly disappointed. ``In my experience, one can translate
three pages even on a very busy day'', he said. And he immediately started to
apologize saying that he did not want to apply pressure, let it be a month, it's OK!
But I perceived it as a challenge. I told to myself: ``24 pages means 8 days.''
I did work hard, and eight days later the paper was translated. ``Good job!'' said Misha, and that was one of the most cherished praises in my entire life.

I will give just one more example. Misha was the chief editor of the Russian translation
of the famous four-volume {\em Linear Partial Differential Operators}\/ by
Lars H\"ormander -- more than two thousand pages of a dense mathematical text.
This is what H\"ormander wrote in his Preface to the Russian edition (my
translation back from Russian):
\begin{quote}
The Russian edition contains several changes, as compared to the English one.
I corrected the errors that I had discovered by myself, while very attentive translators, headed by Professor M.~A.~Shubin made many important clarifications. I was greatly impressed by the meticulous work done by the translating team and I am grateful
for their collaboration -- thanks to it we will be able to introduce improvements
into further English editions.
\end{quote}
So, I was not the only one being impressed by Shubin's style and quality of work, it also
impressed Lars H\"ormander! What more can I say?

\end{input} 

\begin{flushright}
\emph{Maxim Braverman, Arnold Dikansky, Leonid Friedlander, Misha Gromov, Victor Ivrii, Yuri Kordyukov, Peter Kuchment, Vladimir Maz'ya, Robert Mc Owed, Toshikazu Sunada, Aleander Zvonkin}
\end{flushright}

\end{document}